\documentclass[letterpaper, 10 pt, conference]{ieeeconf}  

\IEEEoverridecommandlockouts                              
\usepackage{graphicx}
\usepackage{subcaption}

\usepackage{amsmath,amsfonts,amssymb,amsthm, bm}
\usepackage{hyperref}

\title{\LARGE \bf
Bounding the State Covariance Matrix for a Randomly Switching Linear System with Noise
}

\author{Yongeun Yoon, Corbin Klett, and Eric Feron
\thanks{Y. Yoon is a PhD Student at Georgia Tech.
        {\tt\small yyoon73@gatech.edu}}%
\thanks{C. Klett is a PhD Student at Georgia Tech.
        {\tt\small corbin@ gatech.edu}}%
\thanks{E. Feron is a Professor of Aerospace Engineering at Georgia Tech.
        {\tt\small eric.feron@aerospace.gatech.edu}}%
}

\begin{document}

\maketitle
\thispagestyle{empty}
\pagestyle{empty}

\begin{abstract}

The propagation of a state vector is governed by a set of time-invariant state transition matrices that switch arbitrarily between two values. The evolution of the state is also perturbed by white Gaussian noise with a variance that switches randomly with the state transition relation. The behavior of this system can be characterized by the covariance matrix of the state vector, which is time varying. However, we can bound the set of covariances by comparing the switching system to an augmented system derived with Kronecker algebra. We formulate a matrix optimization problem to compute an ellipsoid that bounds the covariance dynamics, which in turn bounds the state covariance of the set of switching systems subject to white noise. In developing this approach, an invariant ellipsoid for a linear switching affine system is computed along the way.

\end{abstract}

\section{INTRODUCTION}

Methods that can compute guaranteed bounds on the behavior of safety-critical systems are an important part of the control system designer's toolbox. Switching systems are of particular interest since they often appear in practice either intentionally, such as in a hybrid system, or unintentionally, such as when a time delay in the feedback loop effectively changes a system's transition relation. Switching dynamics can also be an abstraction for a system with piecewise-nonlinearities. Many studies have sought to address the need to bound the behavior of switching systems from various perspectives. If the switching strategy is the designer's choice, it may be designed in such a way that a Lyapunov function exists to guarantee stability (\cite{zhang2008robust,kermani2014stability,zhang2009exponential}). Otherwise, the switching order may be constructed by weighted infinite walk \cite{kundu2017stability}. In this research, the switching order is arbitrary. This is motivated by the fact that a system with arbitrarily switching dynamics that is also subject to Gaussian noise can be a template for or a superset of many types of systems for which stability and performance guarantees are difficult to prove, such as time-varying systems, systems with random delays in the feedback loop, and systems subject to various other uncertainties. A relevant study is about the determination of an exponential stability condition for switching systems with an arbitrary switching order \cite{liu2016stability,zhang2017matrix}. However, bounding the trajectories of a switching system involves the determination of an ellipsoid that is guaranteed to be an invariant set. An ellipsoid function can resemble a Lyapunov function since the latter is often chosen to be in a quadratic form. Unlike a Lyapunov function, invariant ellipsoids do not necessarily aim at proving system stability, but, instead, boundedness of the the trajectories followed by the system.

\subsection{Randomly Switching Affine Systems}

First, we study a set of two switching affine systems. In other words, the state propagation is subject to two sets of dynamics, each of which contains a different equilibrium point. A stability guarantee or, at the very least, knowledge of an invariant set for such a system constitutes a preliminary development of this paper. To address these challenges, an invariant set is constructed based on the formulation of a matrix inequality problem in order to provide a guarantee on the state bounds.

\subsection{Random Switching System with Noise}

The randomly switching affine system described above is of the same form of a randomly switching system that propagates a covariance matrix (as the system state) for a randomly switching system with noise. One study seeks to control the covariance matrix for a stochastic discrete-time linear time-varying (LTV) system, steering it from an   initial probabilistic distribution to a desired one \cite{okamoto2018optimal}. A discrete LTV system can be considered to be a switching system with an infinite number of linear affine systems, and steering the covariance matrix from an initial to a desired final condition is a step-wise solution to the problem that this work is interested in. In this research we find the answer to the guaranteed bounding of the covariance for the randomly switching affine system with random affine term for all the steps.

\subsection{Organization} 
First, a formulation for bounding the attractor set of two randomly switching affine systems is developed. To make the problem valid in terms of matrix inequalities, we assume a constant affine term for each system and obtain an ellipsoid that contains the attractor set, that is, all trajectories of the arbitrarily switching affine system as time tends to infinity. This formulation is then generalized to the propagation of a covariance matrix, which can provide a stability guarantee for the switching system with noise. Based on Kronecker algebra, we transform the covariance matrix to a state space representation of switching systems, and apply the same matrix inequality as we did to find the bound of the randomly switching affine systems. Numerical examples demonstrate how the ellipsoid determined through the matrix inequality formulation solves the guaranteed bound problem.

\section{Randomly Switching Affine System}

A simple affine system with dynamics that switch arbitrarily is described by the transition relation
\begin{equation}
\label{eq:rss}
\textbf{x}(k+1) =
\begin{cases}
\textbf{A}_1\textbf{x}(k)+ \textbf{w}_{1} \\
\qquad\text{or} \\
\textbf{A}_2\textbf{x}(k)+ \textbf{w}_{2}
\end{cases} 
\end{equation}
where the real parts of the eigenvalues of $\textbf{A}_i$ all have a magnitude less than one and $\textbf{w}_i$, $i=1,2$ satisfy $\textbf{w}_i=\textbf{A}_i\textbf{x}_{eq,i}-\textbf{x}_{eq,i}$, where $\textbf{x}_{eq,i}$ is an equilibrium point.  System~\ref{eq:rss} does not have an equilibrium in general, unless $\textbf{x}_{eq,1}= \textbf{x}_{eq,2}$. A related notion of equilibrium is that of attractor set~\cite{attractor}, which is the set that the states of system~(\ref{eq:rss}) tend to "converge to". Defining the distance measure between a point $x$ and a set ${\cal S}$ as
\[
d(x,{\cal S}) = \min \left\|y-x\right\|_{y \in {\cal S}},
\]
the attractor ${\cal A}$ of System~(\ref{eq:rss}) is defined as the smallest set such that 
\[
\lim_{k \rightarrow \infty} d(x(k),{\cal A}) = 0.
\]
Little is known in general about the shape or the size of ${\cal A}$. It often features a fractal structure, see Fig.~\ref{fig:rss1} and below discussion for an example.
When it is bounded, we are interested in 
 computing outer approximations in the form of ellipsoidal invariant sets.

\subsection{Computing the Invariant Set}

We search for an invariant set in the form of an ellipse ${\cal E}_{\textbf{P},\textbf{x}_c}$ defined as 
\begin{equation}
{\cal E}_{\textbf{P},\textbf{x}_c} = \left\{\textbf{x} \in {\bf R}^n \; | \; (\textbf{x}-\textbf{x}_c)^T\textbf{P}(\textbf{x}-\textbf{x}_c) \leq 1\right\}
\end{equation}
where $\textbf{x}_c$ is the center of the ellipse and the matrix $\textbf{P}=\textbf{P}^T>0$. In general, an engineer may try to prove asymptotic stability by using $V(\textbf{x}) = (\textbf{x}-\textbf{x}_c)^T\textbf{P}(\textbf{x}-\textbf{x}_c) $ as a candidate Lyapunov function and by verifying the stability condition that $V(\textbf{x}(k+1)) - V(\textbf{x}(k)) \leq 0$, where $\textbf{x}(k+1)$ and $\textbf{x}(k)$ satisfy Equation \ref{eq:rss}. In \eqref{eq:rss}, however, states can propagate in such a way that $V(\textbf{x}(k+1))$ is larger than $V(\textbf{x}(k))$, necessitating an invariant set computation rather than a proof of asymptotic stability; the state will still stay within the ellipsoid if $V(\textbf{x}(k+1))\le 1$ and the absolute stability condition $V(\textbf{x}(k+1))-V(\textbf{x}(k))\le 0$ is not needed. Therefore we consider the following condition for set invariance:
\begin{equation}
\text{if  } V(\textbf{x}(k+1))-V(\textbf{x}(k)) \ge 0,\text{  then  } V(\textbf{x}(k+1)) \le 1
\end{equation}

By the S-Procedure~\cite{boyd1994linear}, this statement is equivalent to 
\begin{equation}
1-V(\textbf{x}(k+1)) - \lambda(V(\textbf{x}(k+1))-V(\textbf{x}(k))) \ge 0
\end{equation}
for some $\lambda>0$. Formulating this statement as a quadratic inequality in the vector $\begin{bmatrix}\textbf{x} & 1\end{bmatrix}^T$, we obtain
\begin{equation}
\label{eq:LMI1}
\!\begin{aligned}
&
\left[\begin{matrix}
(1+\lambda)\textbf{A}_i^T\textbf{P}\textbf{A}_i-\lambda\textbf{P}\\
(1+\lambda)(\textbf{w}_i-\textbf{x}_c)^T\textbf{P}\textbf{A}_i 
\end{matrix}\right.\\
&\quad
\left.\begin{matrix}
(1+\lambda)\textbf{A}_i^T\textbf{P}(\textbf{w}_i-\textbf{x}_c)+\lambda\textbf{P}\textbf{x}_c\\
(1+\lambda)(\textbf{w}_i-\textbf{x}_c)^T\textbf{P}(\textbf{w}_i-\textbf{x}_c)-\lambda\textbf{x}_c^T\textbf{P}\textbf{x}_c-1
\end{matrix}\right]
\end{aligned}\leq0
\end{equation}

A convex program can be solved to find a matrix $\textbf{P}=\textbf{P}^T>0$ which is constrained by a set of $2$ matrix inequalities of the form expressed in \eqref{eq:LMI1} for the randomly switching system~(\ref{eq:rss}). 
By incorporating Schur's lemma to include the initial condition, we derive an additional constraint~\cite{khlebnikov2011}
\begin{equation}
\label{eq:schur}
\begin{bmatrix}
1 & (\textbf{x}_0-\textbf{x}_c)^T \\
\textbf{x}_0-\textbf{x}_c & \textbf{P}
\end{bmatrix} \ge 0
\end{equation}
An SDP solver such as CVX \cite{cvx},\cite{gb08} can find a $\textbf{P}$ satisfying  \eqref{eq:LMI1} and \eqref{eq:schur}. The solution is parameterized by $\lambda$ and $\textbf{x}_0$. By selecting $\textbf{x}_0$ at the centroid of the set of explored states by means of simulations (such a choice is clearly arbitrary and can be changed) and by conducting a line search over $\lambda$, CVX can compute a $\textbf{P}$ of maximum trace (thus an ellipsoid of 
"minimum size") that satisfies the foregoing system of matrix inequalities. 

\subsection{Numerical Example}

A Semi-Definite Programming (SDP) solver such as CVX \cite{cvx,gb08} can find a $\textbf{P}$ satisfying~(\ref{eq:LMI1}) and~(\ref{eq:schur}). The solution is parameterized by $\lambda$ and $\textbf{x}_0$, and therefore with appropriate $\lambda$ the solution $\textbf{P}$ is found. Consider the planar shifted rotating system whose dynamics are defined by~(\ref{eq:rss}). This example is inspired by the fact repeating rotation with a nonzero equilibrium forms a fractal-like pattern~(\cite{peitgen2013fractals}).
\begin{equation}
\begin{aligned}
r=0.9, \quad\theta_1&=0.2,\quad\theta_2=0.1, \\
\textbf{A}_i = r&\begin{bmatrix}
\cos(\theta_i) & -\sin(\theta_i) \\
\sin(\theta_i) & \cos(\theta_i)
\end{bmatrix},\\
\textbf{x}_{eq,1} = \begin{bmatrix}
1 \\
0
\end{bmatrix}&,\quad\textbf{x}_{eq,2} = \begin{bmatrix}
-1 \\
0
\end{bmatrix}, \\
\textbf{w}_{1} = \begin{bmatrix}
-0.118 \\
0.189
\end{bmatrix}, &\quad\textbf{w}_{2} = \begin{bmatrix}
0.104 \\
-0.0899
\end{bmatrix}
\end{aligned} \label{eq:planarrotation}
\end{equation}

\begin{figure}
	\centering
	\includegraphics[width = 0.5\textwidth]{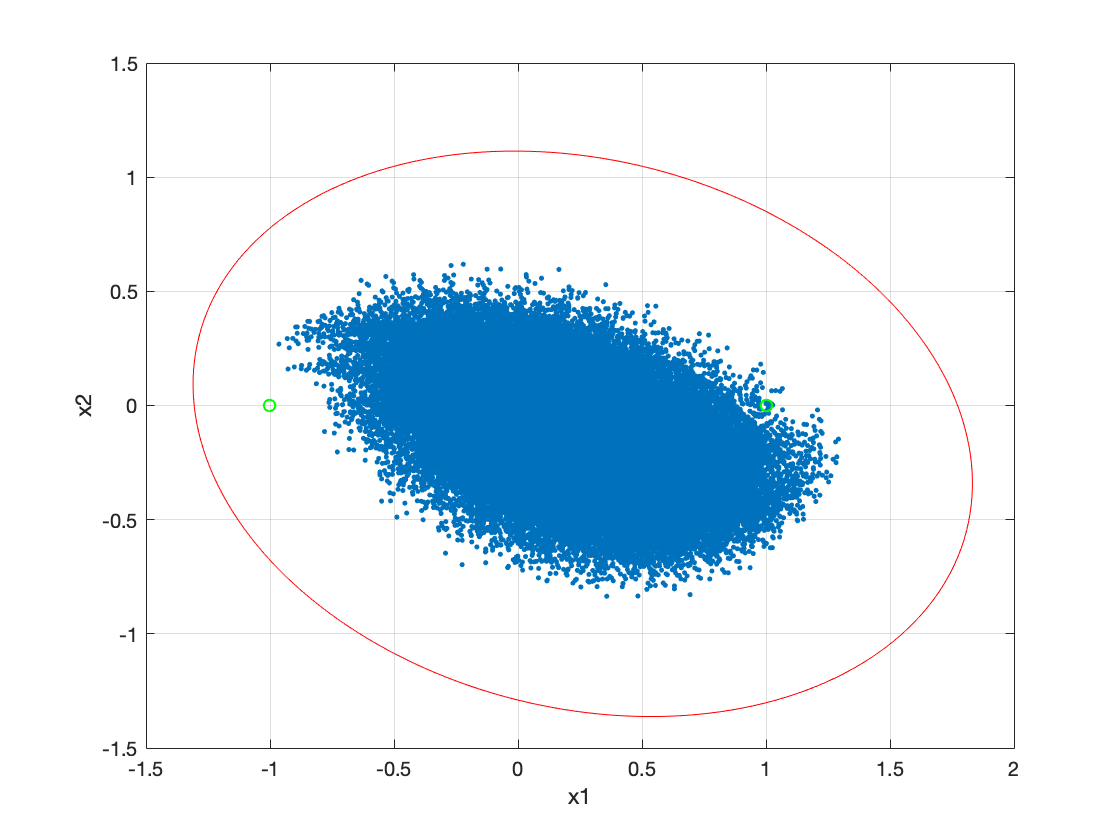}
	\caption{Simulation of system \eqref{eq:rss} with initial condition $\textbf{x}_0=[0 \;0]^T$ and ellipse center $\textbf{x}_c=[0.261 \;-0.124]^T$. The state trajectories form the cluster of points (the \textit{attractor set}) and the invariant ellipse appears around it. Two small circles mark the two equilibria.}
	\label{fig:rss1}
\end{figure}
\noindent
A simulation of 100,000 transitions of the system  \eqref{eq:rss} using the numerical values from \eqref{eq:planarrotation} is shown in Figure \ref{fig:rss1}, as is the invariant ellipse computed using \eqref{eq:LMI1} and \eqref{eq:schur}.


\section{Randomly Switching System with additive Noise}

Consider system similar to \eqref{eq:rss} but with $\textbf{w}_i=\textbf{w}_{i}(k)$, where $\textbf{w}_{i}(k)$, $i=1,2$ is white Gaussian noise with a mean of zero and covariance matrix $\boldsymbol{\Sigma}_i \in \mathbb{R}^{n \text{x} n}$ for a system with state dimension $n$.  An example system that switches between two sets of dynamics is
\begin{equation}
\label{eq:rssn1}
\textbf{x}(k+1) =
\begin{cases}
\textbf{A}_1\textbf{x}(k)+ \textbf{w}_{1}(k) \\ \mbox{or}\\
\textbf{A}_2\textbf{x}(k)+ \textbf{w}_{2}(k) \\
\end{cases} 
\end{equation}
As described in \cite{simon2006optimal}, the propagation of covariance matrices for a discrete transition system with noise $\textbf{w}(k)$ is governed by a discrete-time Lyapunov recursion that converges to a steady state solution $\textbf{P}$ such that $V(\textbf{x})=\textbf{x}^T\textbf{P}\textbf{x}$ is a Lyapunov function for the transition system with no noise. We define the covariance matrix $\textbf{P}_k = \mathbb{E}[\textbf{x}_i(k)\textbf{x}_i(k)^T]$ and a noise covariance matrix as $\textbf{Q}_i=\mathbb{E}[\textbf{w}_{i}(k)\textbf{w}_{i}(k)^T]$. Assuming that the state and the noise are not correlated with each other, 
\begin{equation*}
\begin{aligned}
\textbf{P}_{k+1} &= \mathbb{E}[\textbf{x}_i(k+1)\textbf{x}_i^T(k+1)] \\
&= \mathbb{E}[(\textbf{A}_i\textbf{x}_i(k) + \textbf{w}_{i}(k))(\textbf{A}_i\textbf{x}_i(k) + \textbf{w}_{i}(k))^T] \\
&= \textbf{A}_i\mathbb{E}[\textbf{x}_i(k)\textbf{x}_i^T(k)]\textbf{A}_i^T + \mathbb{E}[\textbf{w}_{i}(k)\textbf{w}_{i}(k)^T]
\end{aligned}
\end{equation*}
resulting in the Lyapunov recursion
\begin{equation}
\textbf{P}_{k+1} = \textbf{A}_i \textbf{P}_k \textbf{A}_i^T + \textbf{Q}_i,	
\label{eq:covpropagate}
\end{equation}
where $i$ is either 1 or 2 at each iteration $k$.
\subsection{Computing a Bound on the State Covariance}
Using Kronecker algebra, the problem of bounding the covariance matrices for the switching system~(ref{eq:rssn1}) with white noise is the same as the problem of computing a bound for the randomly switching affine system~(\ref{eq:rss}). We augment system \eqref{eq:rssn1} by writing\eqref{eq:covpropagate} into the switching affine system
\begin{equation}
\label{eq:kronsys}
\mathcal{P}_{k+1} = \mathcal{A}_i\mathcal{P}_k + \mathcal{Q}_i, \;\; i = 1, 2.
\end{equation}

In the general case, $\mathcal{P}_k=\text{vec}(\textbf{P}_k)$ is the state of system~\ref{eq:kronsys}), $\mathcal{A}_i=\textbf{A}_i\otimes\textbf{A}_i$ the state transition matrix, and $\mathcal{Q}_i=\text{vec}(\textbf{Q}_i)$ is the affine term. The operator $\otimes$ denotes the Kronecker product. Applying \eqref{eq:LMI1} to the augmented sysstem above yields a matrix inequality parameterized by the ellipsoid center $\mathcal{P}_c$ and $\lambda$. If there exists a positive-semidefinite matrix $\mathcal{T}$ for which this matrix inequality holds, then the level set $\mathcal{P}^T\mathcal{T}\mathcal{P}\leq1$ is an invariant set for the sequence of covariance matrices $\textbf{P}_k$, $k=1, \ldots$ of \eqref{eq:rssn1}. The matrix inequality is:
\begin{equation}
\label{eq:LMI2}
\!\begin{aligned}
&
\left[\begin{matrix}
(1+\lambda)\mathcal{A}_i^T\mathcal{T}\mathcal{A}_i-\lambda\mathcal{T}\\
(1+\lambda)(\mathcal{Q}_i-\mathcal{P}_c)^T\mathcal{T}\mathcal{A}_i
\end{matrix}\right.\\
&\quad
\left.\begin{matrix}
(1+\lambda)\mathcal{A}_i^T\mathcal{P}(\mathcal{Q}_i-\mathcal{P}_c)+\lambda\mathcal{T}\mathcal{P}_c\\
(1+\lambda)(\mathcal{Q}_i-\mathcal{P}_c)^T\mathcal{T}(\mathcal{Q}_i-\mathcal{P}_c)-\lambda \mathcal{P}_c^T\mathcal{T}\mathcal{P}_c-1
\end{matrix}\right]
\end{aligned}\leq0
\end{equation}

\subsection{Numerical Example One}

\begin{figure}
	\centering
	\includegraphics[width = .5\textwidth]{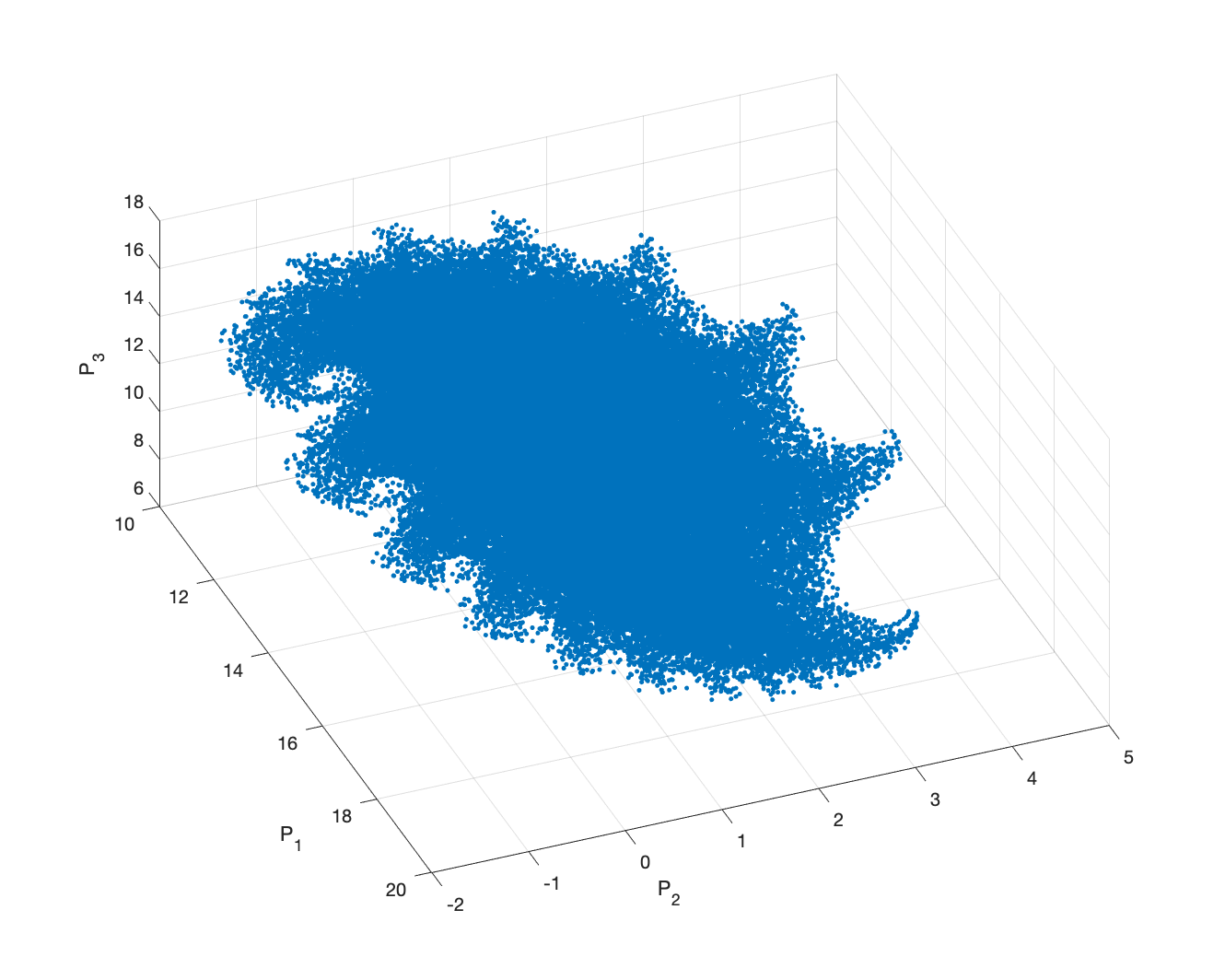}
	\caption{Attractor set of rotating system depicted by simulating 100,000 state transitions }
	\label{fig:fractal1}
\end{figure}

A planar, shifted rotating system of the form \eqref{eq:rssn1} with the state transition matrices described in \eqref{eq:planarrotation} has noise covariance matrices 
\begin{equation}
\label{eq:noise}
\begin{aligned}
\textbf{Q}_1&=
\begin{bmatrix}
2 & 0 \\ 0 & 3
\end{bmatrix},\quad
\textbf{Q}_2&=
\begin{bmatrix}
4 & 0 \\ 0 & 1
\end{bmatrix}
\end{aligned}
\end{equation}
We seek to bound the covariance of the stat. A 100,000 step simulation of system \eqref{eq:rssn1} with an equilibrium at the origin is shown in Figure \ref{fig:fractal1}. The engineer may be interested in generating a bound on the state with a given confidence level, which is guaranteed by the bound on the covariances computed by solving \eqref{eq:LMI2}.

Given that $\textbf{P}$ and the $\textbf{Q}_i$'s are symmetric, we can reduce the dimension of the system \eqref{eq:kronsys} to the dimension of the 4x4 symmetric matrix, which is 3, by eliminating redundant entries produced by symmetry in $\textbf{P}$ and $\textbf{Q}$ and rewriting each $\mathcal{A}_i$ using the entries $a_{jk}$ of $\textbf{A}_i$ as
\begin{equation}
\begin{bmatrix}
a_{11}^2 & 2a_{11}a_{12} & a_{12}^2 \\ a_{21}a_{11} & a_{21}a_{12} + a_{11}a_{22} & a_{22}a_{12} \\ a_{21}^2 & 2a_{21}a_{22} & a_{22}^2
\end{bmatrix}
\end{equation}
If the system $\textbf{x}(k+1)=\textbf{A}\textbf{x}(k)+\textbf{w}(k)$ did not switch dynamics, the propagation of covariance matrices converges to an equilibrium covariance matrix $\textbf{P}_{eq}$ that satisfies the Lyapunov stability criteria $\textbf{P}_{eq}=\textbf{P}_{eq}^T>0$ and $\textbf{A}_i\textbf{P}_{eq}\textbf{A}_i^T-\textbf{P}_{eq}+\textbf{Q}<0$ for some $\textbf{Q}=\textbf{Q}^T>0$. The condition on $\textbf{Q}$ is satisfied by any covariance matrix. Using \eqref{eq:LMI2}, a positive-definite matrix $\mathcal{T}$ is found from which we can construct a quadratic Lyapunov function $V(\mathcal{P})=\mathcal{P}^T\mathcal{T}\mathcal{P}$ that bounds the vector $\mathcal{P}$ for $V(\mathcal{P})\leq 1$. When solving \eqref{eq:LMI2}, the condition number of $\mathcal{T}$ is minimized by placing the center of the ellipse $\mathcal{P}_c$ at the centroid of the dynamics of the set of swtiching systems, which can be done by computing the equilibrium points as $\mathcal{P}_{eq,i}=(\mathcal{A}_i-I)^{-1}\mathcal{Q}_i$ and averaging the two vectors by a weighted factor determined by the inverse of the real parts of the eigenvalues of each $\mathcal{A}_i$.

\begin{figure}
	\centering
	\includegraphics[width = .5\textwidth]{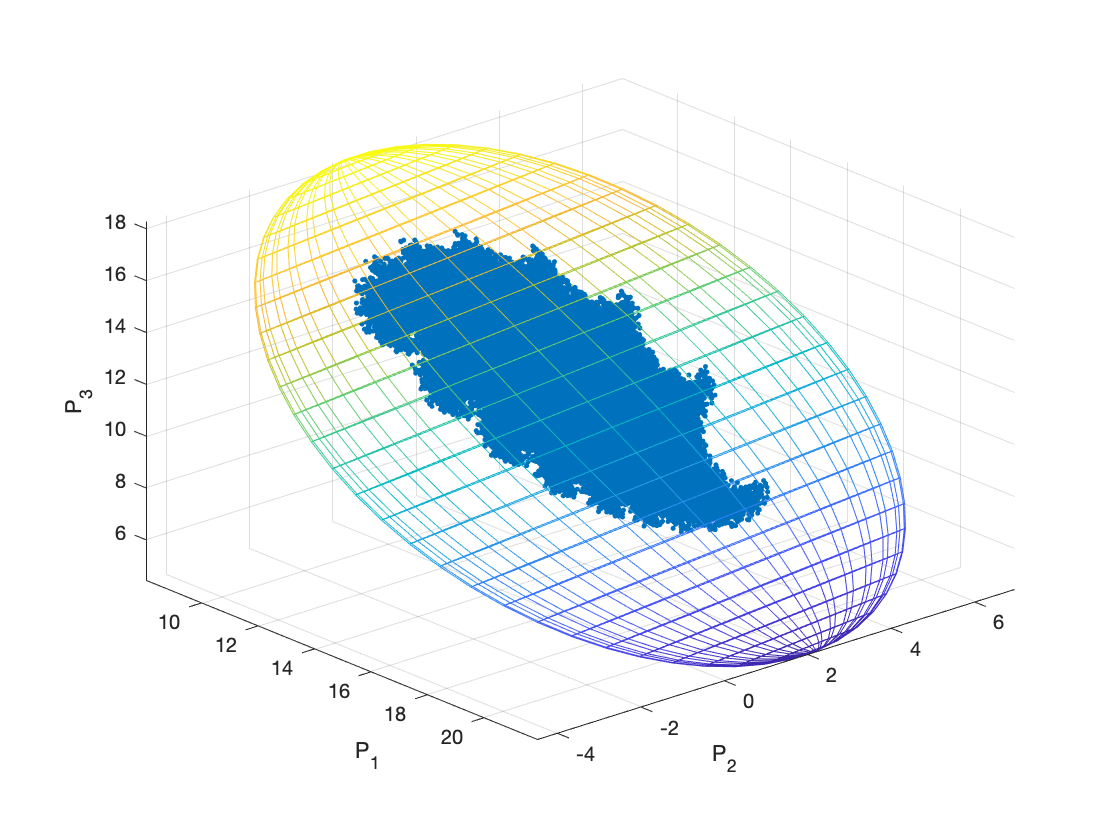}
	\caption{Invariant ellipsoid $\mathcal{P}^T\mathcal{T}\mathcal{P}\leq1$.}
	\label{fig:fractal1_ellipse}
\end{figure}

Finally, a line search on the parameter $\lambda$ is performed in order to maximize the trace of ${\cal T}$. The number of matrix inequalities in the form of \eqref{eq:LMI2} that constrain $\mathcal{T}$ is the same as the number of switching systems. In this numerical example, we compute $\mathcal{T}$ given two matrix inequality constraints. Simulating 100,000 state transitions generates an attractor set that resembles a fractal image, as shown in figure \ref{fig:fractal1}. In this example, all trajectories exist on a plane, but such a reduction in dimensionality does not occur generally (see the next example). The ellipsoid $\mathcal{P}^T\mathcal{T}\mathcal{P} \leq 1$ is shown in figure \ref{fig:fractal1_ellipse}. The ellipsoid squeezes tightly around the plane on which the trajectories lie and contains all possible covariance matrices for the state of system \eqref{eq:rssn1}. The matrix $\mathcal{T}$ is computed in less than a second, with parameters $\lambda=5.21$ and $\mathcal{P}_c=\begin{bmatrix}15.0 & 1.24 & 11.3\end{bmatrix}$ and is found to be
\begin{equation*}
\mathcal{T} = 
\begin{bmatrix}
3.46 & 0.00630 & 3.45 \\ 0.00630 & 0.0306 & 0.00915 \\ 3.45 & 0.00915 & 3.46
\end{bmatrix}
\end{equation*}

\subsection{Numerical Example Two} 

\begin{figure}
	\centering
	\includegraphics[width = .5\textwidth]{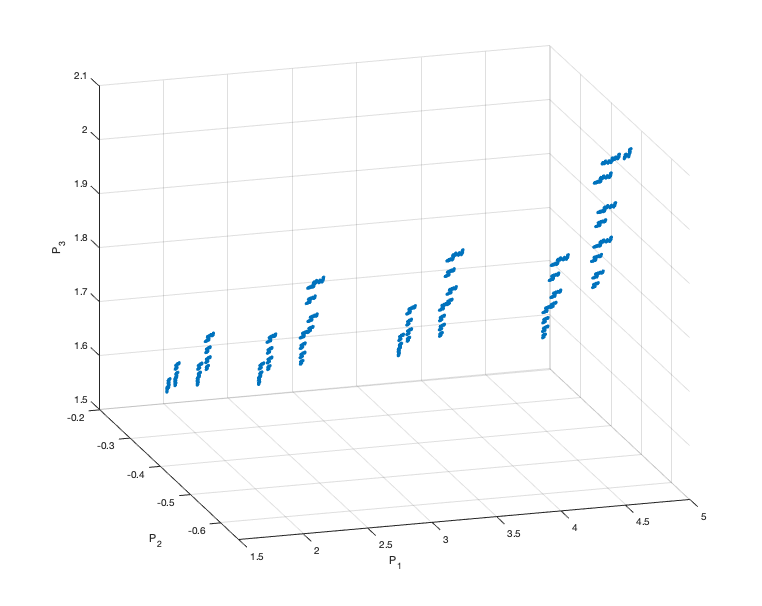}
	\caption{Attractor set of system \eqref{eq:rssn1}, \eqref{eq:erratic}.}
	\label{fig:sys2}
\end{figure}

\begin{figure}
	\centering
	\includegraphics[width = .5\textwidth]{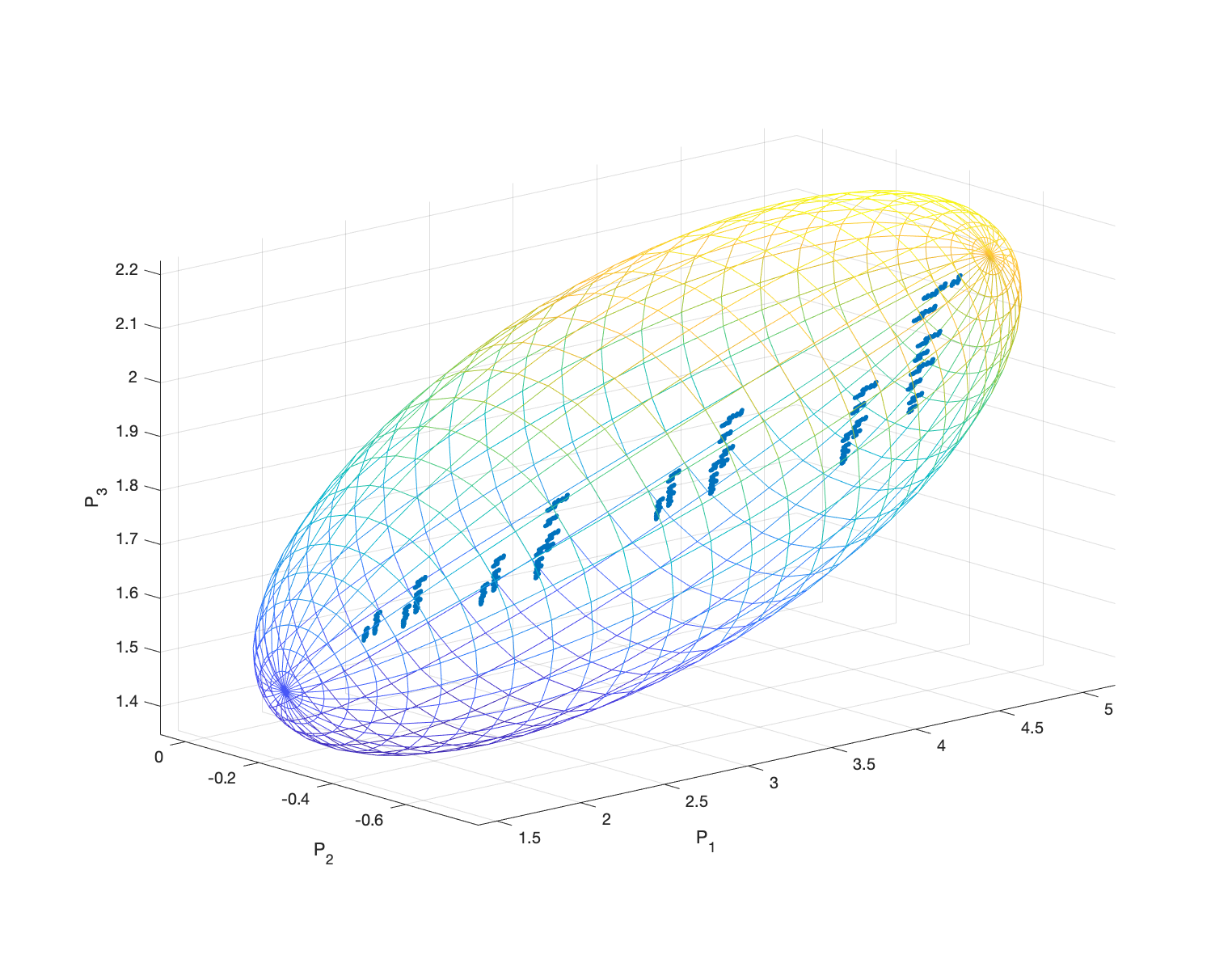}
	\caption{Invariant ellipsoid $\mathcal{P}^T\mathcal{T}\mathcal{P}\leq1$.}
	\label{fig:sys2_ellipse}
\end{figure}

We illustrate the potential of this method by introducing a set of switching systems where the covariance matrices propagate in a far less clustered manner than the previous example (the attractor set does not appear to be contiguous). Consider system \eqref{eq:rssn1} with the dynamics and noise covariances
\begin{equation}
\label{eq:erratic}
\begin{aligned}
\textbf{A}_1 = \begin{bmatrix}
0.7 & -0.7\\
0.2 & 0.7
\end{bmatrix}&, \quad
\textbf{A}_2 = \begin{bmatrix}
0.6 & -0.3\\
0.1 & 0.6
\end{bmatrix} \\
\textbf{Q}_1=
\begin{bmatrix}
1 & 0 \\ 0 & 1
\end{bmatrix}&, \quad
\textbf{Q}_2=
\begin{bmatrix}
1 & 0 \\ 0 & 1
\end{bmatrix}
\end{aligned}
\end{equation}
The simulation of 100,000 transitions is shown in figure \ref{fig:sys2} and the ellipsoid computed by solving the SDP is overlayed in figure \ref{fig:sys2_ellipse}. In this case, the set of possible covariance matrices does not live on a plane.

%

\section{Conclusion}

A problem formulation is presented whereby the attractor of a switching system, which consist of a finite number of linear affine systems, are bounded by an invariant ellipsoid. The bounding condition is transformed into a matrix inequality, and an invariant set can be computed by an SDP solver. Depending on the engineering problem, it even appears that such a set can fit sufficiently tightly around the attractor set produced through simulation. This approach can be used to characterize a randomly switching system with Gaussian noise by constructing an augmented system where the state is the covariance of the original system's state. This research is capable of producing similar results for systems that switch between more than two sets of dynamics by simply increasing the number of inequality constraints in the matrix inequality problem. Issues may arise as the order of the system and number of inequality constraints increases; this is one topic of future study.

\bibliographystyle{ieeetr}
\bibliography{bibliography_yy}

\end{document}